\documentclass{article}

\usepackage[utf8]{inputenc}
\usepackage[T1]{fontenc}

\usepackage{latexsym, amssymb, amsmath, amsthm, multicol,mathrsfs}
\usepackage{graphicx,enumitem}
\usepackage{fullpage}

\usepackage[pdfborder={0 0 0}]{hyperref}
\usepackage{esint}
\usepackage{dsfont}
\usepackage{xcolor}
\usepackage{overpic}
\usepackage{soul}

\newcommand{\cA}{\mathcal{A}}
\newcommand{\cB}{\mathcal{B}}
\newcommand{\cL}{\mathcal{L}}
\newcommand{\cR}{\mathcal{R}}

\newcommand{\N}{\mathbb{N}}

\newcommand{\R}{\mathbb{R}}
\newcommand{\C}{\mathbb{C}}

\renewcommand{\v}{\mathbf{v}}

\renewcommand{\Re}{\mathrm{Re}}

\newcommand{\1}{\mathbf{1}}

\DeclareMathOperator{\Span}{Span}

\DeclareMathOperator*{\esssup}{ess\,sup}
\DeclareMathOperator*{\essinf}{ess\,inf}

\DeclareMathOperator{\suppess}{supp_{ess}}
\DeclareMathOperator{\Div}{div}

\newtheorem{thm}{Theorem}[section]

\newtheorem{lem}[thm]{Lemma}
\newtheorem{cor}[thm]{Corollary}

\author{Bertrand Cloez\thanks{MISTEA, INRAE, Institut Agro, Université de Montpellier, 2 place Pierre Viala, 34060 Montpellier, France.
Email: bertrand.cloez@inrae.fr}
\hspace{1cm} Adil El Abdouni\footnote{Université Paris-Saclay, UVSQ, CNRS, Laboratoire de Mathématiques de Versailles, 78000 Versailles, France. Email: adil.elabdouni76@gmail.com}
\hspace{1cm} Pierre Gabriel\footnote{Institut Denis Poisson, Université de Tours, Université d’Orléans, CNRS, Tours, France. Corresponding author. Email: pierre.gabriel@univ-tours.fr}}

\title{The principal eigenvalue problem for time-periodic nonlocal equations with drift}

\date{}

\begin{document}

\maketitle

\abstract{In this work, we consider a general time-periodic linear transport equation with integral source term.
We prove the existence of a Floquet principal eigenvalue, namely a real number such that the equation rescaled by this number admits nonnegative periodic solutions.
We also prove the exponential attractiveness of these solutions.
The method relies on general spectral results about positive operators.}

\

{\footnotesize{
\noindent\bf MSC2020:} 35B10, 35B40, 35Q92, 35R09, 47B65

}

\section{Introduction}

We consider the following non-autonomous integral partial differential equation
\begin{equation}\label{eq:main}
\partial_t u(t,x) + \Div_x(u(t,x)\v(t,x)) = \int_{\mathbb{R}^d} u(t,y) q(t,y,x)dy + a(t,x)u(t,x)
\end{equation}
where $t\geq0$ and $x\in\R^d$.
The velocity field $\v:\R\times\R^d\to\R^d$, the scalar field $a:\R\times\R^d\to\R$, and the positive kernel $q:\R\times\R^d\times\R^d\to[0,+\infty)$ are supposed to be $T$-periodic in their first variable, for some fixed $T>0$.
Such an equation appears in the modelling of selection-mutations phenomena in varying environments.
In this framework, $u(t,x)\geq0$ represents the repartition of the phenotypical traits $x$ in a population at time $t$.
The transport term represents the drift of the environment.
It is obtained when considering the population in the frame which follows the moving environment.
The zero order term accounts for the selection process through a fitness $a(t,x)$ which is basically the balance between birth and death rates of the trait~$x$ at time~$t$.
The integral term stands for the mutations of the traits, from $y$ to $x$.
A variant consists in modelling the mutations through a Laplace diffusion operator.
We refer for instance to~\cite{Carrere2020,FigueroaIglesias2018,FigueroaIglesias2021,Lorenzi2015} for such models with diffusion approximation of the mutations in periodic environments.

\

In the present paper, we are interested in the spectral properties of Equation~\eqref{eq:main}, and more precisely in the principal eigenvalue problem.
The spectral analysis of time-periodic linear evolution equations dates back to the works of G. Floquet~\cite{Floquet1883}.
Here we aim at proving the existence of a Floquet principal eigenvalue, namely a real value $\lambda_F$ associated to a positive $T$-periodic family ${(f_t)}_{t\in\R}$ such that
\[u(t,x)=e^{\lambda_Ft}f_t(x)\]
satisfies Equation~\eqref{eq:main}.
Our second goal is to investigate the (possibly quantified) exponential stability of this particular solution.

For periodic parabolic equations, the existence of a principal eigenvalue is proved under general assumptions in the book of P. Hess~\cite{Hess1991}, see also~\cite{Huska2008}.
The problem for integral equations attracted attention only more recently.
Existence results for purely integral equations ({\it i.e.} without drift term) in a bounded domain are obtained in the case of a time periodic zero order term and no time dependence of the integral kernel in~\cite{Hutson2008,Rawal2012,Shen2019,Sun2017}.
As far as we know, the case of integral equations with drift was only addressed in~\cite{ElAbdouni2023} in the one dimensional case $d=1$, see also~\cite{Bansaye2020,Clairambault2009,Michel2005} for the singular case of the renewal equation.
Here we propose a general result which is valid in any dimension and under mild assumptions on the coefficients.

\

We make the following hypotheses:
\begin{itemize}[itemindent=4mm,parsep=1mm]
\item[{\bf(H$\v$)}] $\v\in L^\infty_{\mathrm{loc}}(\R,W^{1,\infty}_{\mathrm{loc}}(\R^d))$ and satisfies
\[\v/(1+|x|)\in L^\infty(\R\times\R^d)\qquad\text{and}\qquad \Div_x\v\in L^\infty(\R\times\R^d).\]
For $s,t\in\R$ and $x\in\R^d$, we denote by $X_{t,s}(x)$ the solution to the characteristic equation
\[\left\{\begin{array}{l}
\partial_t X_{t,s}(x)=\v(t,X_{t,s}(x)),
\vspace{2mm}\\
X_{s,s}(x)=x,
\end{array}\right.\]
and we assume that the time spent by these trajectories in balls is uniformly bounded, namely
\begin{equation}\label{eq:DeltaR}
\forall R>0, \quad \sup_{x \in \R^d}\big| \big\{ t  \geq 0, \ |X_{t,0}(x)|<R \big\} \big|  <\infty.
\end{equation}
\item[{\bf(H$a$)}] $a\in L^\infty_{\mathrm{loc}}(\R\times\R^d)$ and tends to $-\infty$ when $|x|\to+\infty$, uniformly in time, in the sense that
\begin{equation}\label{as:asup}
\lim_{R\to+\infty}\esssup_{t\in\R,|x|>R} a(t,x) = -\infty.
\end{equation}
\item[{\bf(H$q$)}] $q\in C_b(\R,L^\infty(\R^d,L^1(\R^d)))$, so that in particular
\begin{equation}\label{as:qsup}
\hat q:=\sup_t\,\esssup_y\int_{\R^d}q(t,y,x)dx<+\infty,
\end{equation}
and it satisfies the non-concentration assumption
\begin{equation}\label{cond Q}
\forall \epsilon>0, \ \exists \eta >0, \ \forall \lvert E \rvert < \eta, \quad \sup_t\,\esssup_y \displaystyle \int_{E} q(t,y,x)dx < \epsilon,
\end{equation}
and the positivity conditions
\begin{equation}\label{as:qpos}
\exists r_0>0,\ \forall t \in \R, \quad q(t,y,x) > 0 \quad \text{for a.e.}\ (x,y)\ \text{s.t.}\ |x-y|<r_0,
\end{equation}
\begin{equation}\label{as:qpos2}
\exists x_0\in\R^d,\  r_1\in(0,r_0), \  q_0 >0, \ \forall t\in\R,\quad \int_{B(x_0,r_1)}\!q(t,y,x)dy \geq q_0 \quad \text{for a.e.}\ x\in B(x_0,r_0).
\end{equation}
\end{itemize}

\

We will also consider the stronger variant

\medskip

\hspace{-3mm}{\bf(H$q+$)} $q\in C_b(\R,L^\infty(\R^d,L^1(\R^d)))$ satisfies~\eqref{cond Q} and the strong positivity condition
\begin{equation}\label{as:qpos+}
\exists r_0,q_0>0,\ \forall t\in \R , \quad q(t,y,x) \geq q_0 \quad \text{for a.e.}\ (x,y)\ \text{s.t.}\ |x-y|<r_0.
\end{equation}

\

Before stating our main theorem, we introduce the dual backward equation of~\eqref{eq:main}, which reads
\begin{equation}\label{eq:main-dual}
-\partial_t \varphi(t,x) = \v(t,x)\cdot\nabla_x\varphi(t,x) + \int_{\mathbb{R}^d} \varphi(t,y) q(t,x,y)dy + a(t,x)\varphi(t,x).
\end{equation}

\begin{thm}\label{thm:main}
Suppose that {\bf(H$\v$)}-{\bf(H$a$)}-{\bf(H$q$)} are met.
Then there exists a unique $\lambda_F\in\R$ and two unique $T$-periodic functions $t\mapsto f_t\in C(\R,L^1_+(\R^d))$ and $t\mapsto\phi_t\in C(\R, L^\infty_+(\R^d))$ such that $\langle\phi_t,f_t\rangle=\|\phi_0\|_{L^\infty}=1$ for all $t\in\R$,
the function $u(t,x)=e^{\lambda_Ft}f_t(x)$ satisfies~\eqref{eq:main}, and $\varphi(t,x)=e^{-\lambda_Ft}\phi_t(x)$ satisfies~\eqref{eq:main-dual}.

Besides, there exist $\rho>0$ and $C\geq1$ such that for any $s\in\R$ and any $u_s\in L^1(\R^d)$, the solution $u(t,x)$ of~\eqref{eq:main} with initial datum $u(s,\cdot)=u_s$ satisfies
\[\big\|e^{-\lambda_F(t-s)}u(t,\cdot)-\langle \phi_s,u_s\rangle f_t\big\|_{L^1}\leq Ce^{-\rho(t-s)}\big\|u_s-\langle \phi_s,u_s\rangle f_s\big\|_{L^1}\qquad\forall t\geq s.\]
If we replace {\bf(H$q$)} by the stronger condition {\bf(H$q+$)}, then the constants $\rho$ and $C$ can be quantified.
\end{thm}

In the case when the coefficients of Equation~\eqref{eq:main} do not depend on time, {\it i.e.} $\v(t,x)=\v(x)$, $a(t,x)=a(x)$, $q(t,x,y)=q(x,y)$, we obtain the following corollary.

\begin{cor}\label{cor:main}
Suppose that the coefficients are time independent and satisfy {\bf(H$\v$)}-{\bf(H$a$)}-{\bf(H$q$)}.
Then there exists a unique $(\lambda_0,f_0,\phi_0)\in\R\times L^1_+(\R^d)\times L^\infty_+(\R^d)$ such that $\langle\phi_0,f_0\rangle=\|\phi_0\|_{L^\infty}=1$,
the function $u(t,x)=e^{\lambda_0t}f_0(x)$ satisfies~\eqref{eq:main}, and $\varphi(t,x)=e^{-\lambda_0t}\phi_0(x)$ satisfies~\eqref{eq:main-dual}.

Besides, there exist $\rho>0$ and $C\geq1$ such that for any $u_0\in L^1(\R^d)$, the solution $u(t,x)$ of~\eqref{eq:main} with initial datum $u(0,\cdot)=u_0$ satisfies
\[\big\|e^{-\lambda_0t}u(t,\cdot)-\langle \phi_0,u_0\rangle f_0\big\|_{L^1}\leq Ce^{-\rho t}\big\|u_0-\langle \phi_0,u_0\rangle f_0\big\|_{L^1}\qquad\forall t\geq s.\]
If we replace {\bf(H$q$)} by the stronger condition {\bf(H$q+$)}, then the constants $\rho$ and $C$ can be quantified.
\end{cor}

Up to our knowledge, the result in Theorem~\ref{thm:main} is new in the literature, extending the result of~\cite{ElAbdouni2023} to higher dimensions and more general coefficients.
In the case of time independent coefficients, Corollary~\ref{cor:main} extends~\cite[Theorem~4.2]{Velleret2023} to general drift rates, \cite[Theorem~4.1]{Li2017} to the whole space in any dimension, and~\cite[Theorem~2.1]{Cloez2020} to higher dimensions.

\medskip

In~\cite{ElAbdouni2023}, the author uses the non-conservative Harris theory by adapting the results of~\cite{Bansaye2022,Cloez2020} to the periodic setting.
Here we use another approach by first proving the existence of the Floquet principal eigenvalue and the associated periodic eigenfunctions by means of spectral arguments.
To this end, in Section~\ref{sec:operator} we extend the methodology of~\cite{FGM} to derive operator theoretic results that are applicable to periodic semiflows.
These results are then applied to the semiflow generated by Equation~\eqref{eq:main} in Section~\ref{sec:EDP}.
This spectral approach also provides the non-constructive exponential asymptotic stability of Floquet's solutions.
We then apply Harris's method for quantifying the exponential rate of convergence under the strengthened assumption {\bf(H$q+$)}.

\

\section{Some operator theoretic results}\label{sec:operator}

In this section, we prove general results about positive operators in $L^p$ spaces, that are most probably extendable to some other Banach lattices.
The approach is inspired from the material presented and developed in~\cite{FGM,Lions}.
We consider $X=L^p(E,\mathscr E,\mu)$ the Lebesgue space of functions associated to the Borel $\sigma$-algebra~$\mathscr E$, a positive $\sigma$-finite measure $\mu$ and an exponent $p\in[1,\infty)$,
or $X=L^p_m(E)$ its weighted variant with $m:E\to(0,\infty)$ a measurable weight function.

\medskip

We denote by $X_+$ the standard positive cone of $X$ made of the almost everywhere nonnegative functions of $X$.
For such a function $f$ we write $f\geq0$, and we thus have $X_+=\{f\in X,\ f\geq0\}$.
For a function $f\in X$ we set $f_+=\max(f,0)$ and $f_-=\max(-f,0)$, so that $f=f_+-f_-$. 
An operator $U:X\to X$ is said to be positive, and we write $U\geq0$, if
\[f\in X_+\ \quad\implies\quad Uf\in X_+.\]
We will sometimes abuse notations by writing $U:X_+\to X_+$ to mean that $U$ is a positive operator on $X$.
For an operator $U\geq0$ we have
\begin{equation}\label{eq:ineqU>0}
|Uf|\leq U|f|
\end{equation}
for any $f\in X$, where $|f|=f_++f_-$ is the absolute value of $f$.
We use the notation $f>0$ for meaning that $f\in X_+\setminus\{0\}$, and for a function which is strictly positive almost everywhere we write $f\gg0$.
We also introduce the cone $X_{++}=\{f\in X,\ f\gg0\}$ and we say that an operator $U$ on $X$ is strictly positive if
\[f\in X_+\setminus\{0\}\ \quad\implies\quad Uf\in X_{++}.\]
For such an operator we write $U>0$ or, abusing notations, $U:X_+\setminus\{0\}\to X_{++}$.

\medskip

For studying the spectrum of an operator $U$ in $X$, we need to define the complex space $X_\C$ made of the complex valued functions such that the modulus $|f|$ belongs to $X$.
For avoiding confusions we will sometimes write $X_\R$ for the real space $X$.
A function $f\in X_\C$ is uniquely written as $f=g+ih$ with $g,h\in X_\R$,
and we have $|f|=\sqrt{g^2+h^2}$ but also
\[|f|=\sup_{0\leq\theta\leq2\pi}\Re(e^{-i\theta}f)=\sup_{0\leq\theta\leq2\pi}(g\cos\theta+h\sin\theta).\]
From this last identity we readily extend, for a positive operator $U$, the inequality~\eqref{eq:ineqU>0} from $X=X_\R$ to $X_\C$.

\medskip

For $U:X\to X$ a bounded operator, we denote by
\[\rho(U):=\{\lambda\in\C,\ \lambda-U\ \text{is bijective}\}\]
its resolvent set and by $\sigma(U)=\C\setminus\rho(U)$ its spectrum.
We denote by
\[\sigma_p(U):=\{\lambda\in\C,\ \lambda-U\ \text{is not injective}\}\subset\sigma(U)\]
the point spectrum of $U$, namely the set of eigenvalues of $U$,
and by
\[r(U):=\sup\{\lambda\in\sigma(U)\}=\inf\{r>0,\ B^c(0,r)\subset\rho(U)\}\]
the spectral radius of $U$.
We recall the spectral radius formula  (see \cite[Thm.~10.13]{Rudin}  for instance)
\begin{equation}\label{eq:SRF}
r(U)=\lim_{n\to\infty}\|U^n\|^{1/n},
\end{equation}
which ensures in particular that $r(U^*)=r(U)$, for $U^*:X'\to X'$ the dual operator of $U$.
Indeed,
\[\|U^*\|=\sup_{\|\phi\|_{X'}=1}{\|U^*\phi\|}_{X'}=\sup_{{\|f\|}_X=\|\phi\|_{X'}=1}\langle U^*\phi,f\rangle=\sup_{{\|f\|}_X=\|\phi\|_{X'}=1}\langle \phi,Uf\rangle=\sup_{{\|f\|}_{X}=1}{\|Uf\|}_{X}=\|U\|,\]
where $\langle\cdot,\cdot\rangle$ stands for the standard duality bracket $\langle\cdot,\cdot\rangle_{X',X}$, and the forth equality is a consequence of the Hahn-Banach theorem.
We also recall that an operator $U$ on $X$ is said to be power compact if $U^\ell$ is compact for some integer $\ell$.

\medskip

In what follows, we consider $U:X_+\to X_+$ a positive bounded operator on $X$ which satisfies $r(U)>0$.

\begin{thm}\label{thm:U1}
Assume that $U$ admits the splitting
\begin{equation}\label{eq:splittingU}
U=W+K\quad \text{with}\ \|W\|<r(U)\ \text{and}\ K(\lambda-W)^{-1}\ \text{power compact for any}\ |\lambda|>\|W\|.
\vspace{-2.5mm}
\end{equation}
Then
\vspace{-1.5mm}
\begin{enumerate}[label=\roman*),itemsep=-3pt]
\item there exist $f_0\in X_+\setminus\{0\}$ and $\phi_0\in X'_+\setminus\{0\}$ such that
\begin{equation}\label{eq:vpU}
Uf_0=r(U)f_0\quad \text{and}\quad U^*\phi_0=r(U)\phi_0,
\end{equation}
\item for any $\kappa$ such that $\|W\|<\kappa<r(U)$, the set $\sigma(U)\cap B^c(0,\kappa)$ is a finite subset of $\sigma_p(U)$.
\end{enumerate}
\end{thm}

Let us point out that the power compactness assumption of $K(\lambda-W)^{-1}$ in~\eqref{eq:splittingU} corresponds to the $W$-power compactness of $K$ in the terminology introduced by J. Voigt in~\cite{Voigt1980}.
The property~{\it ii)} of Theorem~\ref{thm:U1} is then a direct consequence of~\cite[Corollary~1.4]{Voigt1980}.

\begin{proof}[Proof of Theorem~\ref{thm:U1} i)]
The spectral radius formula~\eqref{eq:SRF} implies that for any real number $\lambda>r(U)$ the series $\sum_{k\geq0}\lambda^{-k}U^k$ is convergent, which ensures that $\lambda-U$ is invertible with inverse given by
\[\cR(\lambda):=(\lambda-U)^{-1}=\lambda^{-1}\sum_{k=0}^\infty\lambda^{-k}U^k.\]
Since $U\geq0$ and $r(U)\geq0$, this inverse is also a positive operator.

\smallskip

Now we prove the classical result that for any $\lambda$ in $\rho(U)$, the resolvent set of $U$, we have
\[\|(\lambda-U)^{-1}\|\geq\frac{1}{d(\lambda,\sigma(U))}.\]
For any $\lambda\in\rho(U)$ and $\mu\in\C$ we have
\[(\mu-U)=\big[I-(\lambda-\mu)\cR(\lambda)\big](\lambda-U).\]
Since if $|\mu-\lambda|<1/\|\cR(\lambda)\|$ the series $\sum(\lambda-\mu)^n\cR(\lambda)^n$ converges to $\big[I-(\lambda-\mu)\cR(\lambda)\big]^{-1}$ we get
\begin{align*}
    (\mu-U)\sum_{n=0}^\infty(\lambda-\mu)^n\cR(\lambda)^{n+1}&
    =\big[I-(\lambda-\mu)\cR(\lambda)\big](\lambda-U)\sum_{n=0}^\infty(\lambda-\mu)^n\cR(\lambda)^{n+1}\\
    &=\big[I-(\lambda-\mu)\cR(\lambda)\big]\sum_{n=0}^\infty(\lambda-\mu)^n\cR(\lambda)^n=I.
\end{align*}
This ensures that $B(\lambda,1/\|\cR(\lambda)\|)\subset \rho(U)$ or equivalently $d(\lambda,\sigma(U))\geq1/\|\cR(\lambda)\|.$

\smallskip

Now let $(\mu_n)\subset \overline B^c(0,r(U))\subset \rho(U)$ be a sequence such that $d(\mu_n,\sigma(U))\to0$, and so $\|\cR(\mu_n)\|\to+\infty$.
Consider the sequence $\lambda_n=|\mu_n|>r(U)$ which satisfies $\lambda_n\to r(U)$.
Since $|Uf|\leq U|f|$, we have that
\[|\cR(\mu_n)f|=|\mu_n|^{-1}\bigg|\sum_{k=0}^\infty\mu_n^{-k}U^kf\bigg|
\leq\lambda_n^{-1}\sum_{k=0}^\infty\lambda_n^{-k}U^k|f|=\cR(\lambda_n)|f|.\]
This ensures that $\|\cR(\mu_n)\|\leq\|\cR(\lambda_n)\|$ and consequently $\|\cR(\lambda_n)\|\to+\infty$.
This means that there exist two sequences $(f_n)$ and $(g_n)$ such that
\[f_n=\cR(\lambda_n)g_n,\quad \|f_n\|\to\infty\quad\text{and}\quad \|g_n\|\leq1.\]
By splitting $g_n = g_n^+-  g_n^-$, we get
\[f_n = \cR(\lambda_n) g_n^+ - \cR(\lambda_n) g_n^-\]
with 
\[\| g_n^\pm \| \le 1 \quad \hbox{and}\quad (\| \cR(\lambda_n) g_n^+ \| \to \infty \text{ or } \| \cR(\lambda_n) g_n^- \| \to \infty).\]
Changing notations, we thus have the existence of two sequences $(f_n)$ and $(g_n)$ such that
\[f_n=\cR(\lambda_n)g_n,\quad f_n\geq0,\quad g_n\geq0,\quad\|f_n\|\to\infty\quad\text{and}\quad \|g_n\|\leq1.\]
We define $h_n:=(\lambda_n-W)f_n$, which also satisfies $\|h_n\|\to\infty$ since
\[\|h_n\|\geq\|(\lambda_n-W)^{-1}\|^{-1}\|f_n\|\geq(\lambda_n-\|W\|)\|f_n\|\geq(r(U)-\|W\|)\|f_n\|\]
with $r(U)>\|W\|$.
Using that $U=W+K$ and setting $\hat h_n := h_n / \| h_n \|$ and $\varepsilon_n := g_n/ \| h_n \|$, we get that
\[\hat h_n = K(\lambda_n-W)^{-1}\hat h_n + \varepsilon_n,\quad \|\hat h_n\|=1\quad\text{and}\quad \|\varepsilon_n\|\to0.\]
Iterating this equality, we also have for any integer $\ell\geq1$
\[\hat h_n = \big[K(\lambda_n-W)^{-1}\big]^\ell\hat h_n+\widetilde\varepsilon_n\]
with $\|\widetilde\varepsilon_n\|\to0$.
Choosing $\ell$ large enough so that $\big[K(r(U)-W)^{-1}\big]^\ell$ is compact,
we deduce the existence of $h_0\in X$ and a subsequence of $(\hat h_n)$, not relabeled, such that
\[\big\|\big[K(r(U)-W)^{-1}\big]^\ell\hat h_n-h_0\big\|\to0.\]
Due to the continuity of $\lambda\mapsto (\lambda-W)^{-1}$ we also have $\|\big[K(\lambda_n-W)^{-1}\big]^\ell\hat h_n-h_0\|\to0$ and consequently $\|\hat h_n-h_0\|\to0$.
In particular $h_0\neq0$, since $\|h_0\|=\lim\|\hat h_n\|=1$, and
\[h_0=K(r(U)-W)^{-1}h_0.\]
Setting $f_0=(r(U)-W)^{-1}h_0$ we thus have
\[Uf_0=r(U)f_0,\quad f_0\neq0\quad \text{and}\quad f_0\geq0,\]
where the last property comes from the fact that $f_0=\lim f_n/\|h_n\|$ with $f_n\geq0$.

\smallskip

Arguing similarly for the dual operator $U^*$ we have the existence of $(\phi_n)$ and $(\psi_n)$ such that
\[\phi_n=\cR(\lambda_n)^*\psi_n,\quad \phi_n\geq0,\quad \psi_n\geq0,\quad\|\phi_n\|\to\infty\quad\text{and}\quad \|\psi_n\|\leq1.\]
Setting $\hat\phi_n=\phi_n/\|\phi_n\|$ and $\varepsilon_n=\psi_n/\|\phi_n\|$ we get that
\[\hat \phi_n = \big[(\lambda_n-W^*)^{-1}K^*\big]^\ell\hat \phi_n+\varepsilon_n,\quad \|\hat \phi_n\|=1\quad\text{and}\quad \|\varepsilon_n\|\to0.\]
For $\ell$ such that $\big[K(\lambda_n-W)^{-1}\big]^\ell$ is compact, its dual $\big[(\lambda_n-W^*)^{-1}K^*\big]^\ell$ is also compact.
We deduce the existence of a subsequence of $(\hat\phi_n)$ which converges to a limit $\phi_0\in X'_+$ satisfying
\[U^*\phi_0=r(U)\phi_0\quad \text{and}\quad \|\phi_0\|=1.\]
\end{proof}

For the next results, we define $\widetilde U:=\frac{1}{r(U)}U$, which thus satisfies $r(\widetilde U)=1$.
We also recall the notation $N(A)=\{f\in X,\ Af=0\}$ for the null space of a linear operator $A$.

\begin{thm}\label{thm:U2}
Assume that $U:X_+\setminus\{0\}\to X_{++}$ and that there exist $f_0\in X_{+}\setminus\{0\}$ and $\phi_0\in X'_{+}\setminus\{0\}$ satisfying~\eqref{eq:vpU}, or equivalently such that $\widetilde Uf_0=f_0$ and $\widetilde U\phi_0=\phi_0$.
Then
\vspace{-1mm}
\begin{enumerate}[label=\roman*),itemsep=0pt]
\item $f_0\in X_{++}$ and $\phi_0\in X'_{++}$,
\item $N(\widetilde U-I)=\Span(f_0)$ and $N(\widetilde U^*-I)=\Span(\phi_0)$,
\item the peripheral point spectrum of $\widetilde U$, namely the set of eigenvalues of $\widetilde U$ with modulus equal to $r(\widetilde U)=1$, is a multiplicative subgroup of $\mathbb S^1:=\{z\in\C,\ |z|=1\}$.
\end{enumerate}
\end{thm}

\begin{proof}[Proof of Theorem~\ref{thm:U2} i)]
The strict positivity of $f_0$ follows immediately from the strict positivity assumption made on $U$.
For $\phi_0$, we observe that for any $f\in X_+\setminus\{0\}$
\[\langle\phi_0,f\rangle=\langle\widetilde U^*\phi_0,f\rangle=\langle\phi_0,\widetilde Uf\rangle>0\]
since $\phi_0\in X_+\setminus\{0\}$ and $\widetilde Uf\in X_{++}$ due again to the strict positivity of $U$.
Since $f$ is chosen arbitrarily in $X_+\setminus\{0\}$, we deduce that $\phi_0$ must be in $X'_{++}$.
\end{proof}

\begin{proof}[Proof of Theorem~\ref{thm:U2} ii)]
Let $f\in X\setminus\{0\}$ such that $\widetilde Uf=f$.
First we observe that, by positivity of $U$,
\[|f|=|\widetilde Uf|\leq\widetilde U|f|.\]
This inequality is actually an equality, otherwise we would have, due to the strict positivity of $\phi_0$,
\[\langle\phi_0,|f|\rangle<\langle\phi_0,\widetilde U|f|\rangle=\langle\widetilde U^*\phi_0,|f|\rangle=\langle\phi_0,|f|\rangle\]
and a contradiction.
We thus have $\widetilde Uf=f$ and $\widetilde U|f|=|f|$, from what we deduce that
\[\widetilde Uf_\pm=f_\pm\]
by writing $f_\pm=\frac12(|f|\pm f)$.
The strict positivity of $U$ implies that $f_\pm\in L^1_{++}$ or $f_\pm=0$, and thus $f_+\in X_{++}$ or $f_-\in X_{++}$, since $f\neq0$.
Without loss of generality we may assume that $f_+\in X_{++}$, meaning that $f_+(x)>0$ for almost everywhere.
This necessarily implies that $f_-(x)=0$ almost everywhere and thus $f=f_+\in X_{++}$.
We now introduce
\[g=\langle\phi_0,f_0\rangle f-\langle\phi_0,f\rangle f_0,\]
which satisfies $\widetilde Ug=g$ and $\langle\phi_0,g\rangle=0$.
This enforces $g=0$.
Indeed, otherwise, arguing as above we would obtain that either $g\in X_{++}$ or $-g\in X_{++}$, both contradicting $\langle\phi_0,g\rangle=0$.
We thus have proved that $f\in\Span(f_0)$ and so $N(\widetilde U-I)=\Span(f_0)$.
The same arguments ensure that $N(\widetilde U^*-I)=\Span(\phi_0)$.
\end{proof}

Before proving the last point of Theorem~\ref{thm:U2}, we introduce some additional notations.
For $f\gg0$ we define
\[X_f:=\big\{g\in X,\ \exists C>0,\ |g|\leq Cf\big\}\]
and for $g\in X_f$ we set
\[{[g]}_f:=\inf\big\{C>0,\ |g|\leq Cf\big\}.\]
The proof of Theorem~\ref{thm:U2}~{\it iii)} is based on the following result, taken from~\cite[Lemma~5.3]{FGM}.

\begin{lem}[\cite{FGM}]\label{lem:Q>0}
Let $f\in X_{++}$ and consider a linear operator $Q:X_f\to X_f$ such that $Qf=f$ and ${[Qg]}_f\leq{[g]}_f$ for any $g\in X_f$. Then $Q\geq0$.
\end{lem}

\begin{proof}[Proof of Theorem~\ref{thm:U2} iii)]
{\it Step 1.}
Consider $f\in X_\C\setminus\{0\}$ such that $\widetilde Uf=e^{i\alpha}f$ for some $\alpha\in(0,2\pi)$.
First we remark that
\[|f|=|e^{i\alpha} f|=|\widetilde Uf|\leq\widetilde U|f|.\]
For the same reason as in the proof of point {\it ii)}, this inequality is actually an equality and we have $\widetilde U|f|=|f|$.
Since $N(\widetilde U-I)=\Span(f_0)$ we deduce that $|f|\in X_{++}$.
Now we define the operator $Q$ by
\[Qg:=\frac{\bar f}{|f|}e^{-i\alpha}\widetilde U\bigg(\frac{fg}{|f|}\bigg).\]
On the one hand, we have
\[Q|f|=\frac{\bar f}{|f|}e^{-i\alpha}\widetilde Uf=|f|=\widetilde U|f|.\]
On the other hand, we have
\[|Qg|\leq \bigg|\widetilde U\bigg(\frac{fg}{|f|}\bigg)\bigg|\leq\widetilde U|g|,\qquad\forall g\in L^1,\]
which, by positivity of $\widetilde U$, yields
\[|Qg|\leq{[g]}_f\widetilde U|f|={[g]}_f|f|,\qquad\forall g\in X_f.\]
From Lemma~\ref{lem:Q>0} applied to $|f|\gg0$, we deduce that $Q\geq0$ on $X_{|f|}$ and then on $X=\overline{X_{|f|}}$.
As a consequence, $0\leq Qg=|Qg|\leq\widetilde Ug$ for any $g\geq0$.
In other words we have $\widetilde U-Q\geq0$, and then $\widetilde U^*-Q^*\geq0$.
We must have $\widetilde U^*-Q^*=0$.
Otherwise, there would exist $\psi\in X'_+\setminus\{0\}$ such that $(\widetilde U^*-Q^*)\psi\in X'_+\setminus\{0\}$,
and we find a contradiction by computing
\[0<\langle(\widetilde U^*-Q^*)\psi,|f|\rangle=\langle\psi,(\widetilde U-Q)|f|\rangle=0.\]
We have established that $\widetilde U=Q$, and so
\[\widetilde U\bigg(\frac{fg}{|f|}\bigg)=e^{i\alpha}\frac{f}{|f|}\widetilde Ug,\qquad \forall g\in X.\]

\smallskip

{\it Step 2.}
Consider now $\alpha,\beta\in(0,2\pi)$ and $f,g\in X\setminus\{0\}$ such that $\widetilde Uf=e^{i\alpha}f$ and $\widetilde Ug=e^{i\beta}g$ and define $h=\dfrac{fg}{|fg|}f_0$.
Using the step 1 we get
\[\widetilde Uh=e^{i\alpha}\frac{f}{|f|}\widetilde U\bigg(\frac{gf_0}{|g|}\bigg)=e^{i(\alpha+\beta)}\frac{fg}{|fg|}\widetilde U f_0=e^{i(\alpha+\beta)}h.\]
This proves that $e^{i\alpha}e^{i\beta}\in\sigma_p(\widetilde U)$, and the proof is complete.
\end{proof}

We now state a theorem which is obtained by combining the results of Theorems~\ref{thm:U1} and~\ref{thm:U2}.

\begin{thm}\label{thm:U3}
Assume that $U:X_+\setminus\{0\}\to X_{++}$ and that $U$ admits the splitting~\eqref{eq:splittingU}.
Then
\begin{enumerate}[label=\roman*)]
\item there exists $(f_0,\phi_0)\in X_{++}\times X'_{++}$ such that $\widetilde Uf_0=f_0$ and $\widetilde U^*\phi_0=\phi_0$,
\item $N(\widetilde U-I)=\Span(f_0)$ and $N(\widetilde U^*-I)=\Span(\phi_0)$,
\item the peripheral point spectrum of $\widetilde U$ is reduced to $\{1\}$,
\item there exists $\zeta\in(0,1)$ and $C\geq1$ such that for all $n\in\N$ and all $f\in {\{\phi_0\}}^\perp:=\big\{g\in L^1,\,\langle\phi_0,g\rangle=0\big\}$
   \begin{equation}\label{eq:discrete-ergo-U}
   \big\|\widetilde U^nf\big\|_X\leq C\, \zeta^n\left\|f\right\|_X.
   \end{equation}
\end{enumerate}
\end{thm}

The points {\it i)} and {\it ii)} are direct consequences of Theorems~\ref{thm:U1} and~\ref{thm:U2}, and we thus only prove the items {\it iii)} and {\it iv)}. 

\begin{proof}[Proof of Theorem~\ref{thm:U3} iii)]
From Theorem~\ref{thm:U1}~{\it ii)} and Theorem~\ref{thm:U2} {\it iii)} we know that the peripheral spectrum of $\widetilde U$ is a finite subgroup of $\mathbb S^1$.
Suppose by contradiction that it is not reduced to the trivial subgroup $\{1\}$, and so that there exists an integer $\ell\geq2$ such that $e^{2i\pi/\ell}\in\sigma_p(\widetilde U)$.
This means that
\begin{equation}\label{eq:complexeigen}
\exists f=g+ih\in X_\C\setminus\{0\},\qquad\widetilde Uf=e^{2i\pi/\ell}f.
\end{equation}
This identity has two incompatible implications.

\smallskip

On the one hand, by iteration, we get that $\widetilde U^\ell f=f$.
Taking the real and imaginary parts, we get that $\widetilde U^\ell g=g$ and $\widetilde U^\ell h=h$.
Since $f\neq0$, we must have $g\neq0$ or $h\neq0$ (and actually both if $\ell\geq3$).
Assume without loss of generality that $g\neq0$.
Then the strict positivity assumption on $U$ implies that $g\gg0$.

\smallskip

On the other hand, testing~\eqref{eq:complexeigen} against $\phi_0$ gives
\[\langle\phi_0,f\rangle=\langle\widetilde U^*\phi_0,f\rangle=\langle\phi_0,\widetilde U f\rangle = e^{2i\pi/\ell}\langle\phi_0,f\rangle.\]
Since $e^{2i\pi/\ell}\neq1$, we get that $\langle\phi_0,f\rangle=0$ and as a consequence $\langle\phi_0,g\rangle=0$.
This is a contradiction with the facts that $g\gg0$ and $\phi_0\gg0$.
\end{proof}

\begin{proof}[Proof of Theorem~\ref{thm:U3} iv)]
Since $\phi_0$ is an eigenvector of $\widetilde U^*$, the subspace $\{\phi_0\}^\perp$ of $X$ is invariant under $\widetilde U$.
Denote by $\widetilde U^\perp$ the restriction of $\widetilde U$ to this invariant subspace.
Since $f_0\gg0$ and $\phi_0\gg0$, we have that $f_0\not\in\{\phi_0\}^\perp$ and the point {\it ii)} then implies that $N(\widetilde U^\perp-I)=\{0\}$, which in turns yields that $1\not\in\sigma(\widetilde U^\perp)$.
Using Theorem~\ref{thm:U1}~{\it ii)} and Theorem~\ref{thm:U3}~{\it iii)}, we infer that $r(\widetilde U^\perp)<1$.
The spectral radius formula~\eqref{eq:SRF} applied to $\widetilde U^\perp$ then ensures that for any $\zeta\in(r(\widetilde U^\perp),1)$ there exists $C\geq1$ such that
\[\big\|\big(\widetilde U^\perp\big)^n\big\|\leq C\zeta^n,\qquad\forall n\in\N.\]
This is exactly~\eqref{eq:discrete-ergo-U}.
\end{proof}

Now we give a corollary of Theorem~\ref{thm:U3} in the $L^1$ framework that we will use for proving Theorem~\ref{thm:main}.

\begin{cor}\label{cor:abstract}
    Let $S:L^1_+\to L^1_+$ be a positive bounded operator and suppose that for some $k\in\N$ and $\kappa_0>0$
    \begin{enumerate}
        \item $\exists g_0\in L^1_+\setminus\{0\}$ such that $Sg_0\geq \kappa_0 g_0$, or $\exists \psi_0\in L^\infty_+\setminus\{0\}$ such that $S^*\psi_0\geq \kappa_0 \psi_0$,
        \item $S^k=W+K$ with $\|W\|<\kappa_0^k$ and $K$ weakly compact in $L^1$,
        \item $\forall f\in L^1_+\setminus\{0\},\ S^kf(x)>0\ \text{for almost every}\ x\in\R^d$.
    \end{enumerate}
   Then there exists a unique triplet $(\Lambda_0,f_0,\phi_0)\in \R\times L^1_+(\R^d)\times L^\infty_+(\R^d)$ such that
   \[Sf_0=\Lambda_0f_0,\quad S^*\phi_0=\Lambda_0\phi_0,\quad \text{with}\quad \|\phi_0\|_{L^\infty}=\langle\phi_0,f_0\rangle=1.\]
   Moreover, $\Lambda_0>0$, $f_0,\phi_0\gg0$, and there exists $\zeta\in(0,1)$ and $C\geq1$ such that for $n\in\N$ and all $f\in L^1(\R^d)$
   \begin{equation}\label{eq:discrete-ergo}
   \left\|\Lambda_0^{-n}S^nf-\langle\phi_0,f\rangle f_0\right\|_{L^1}\leq C \zeta^n\left\|f-\langle\phi_0,f\rangle f_0\right\|_{L^1}.
   \end{equation}
\end{cor}

\begin{proof}[Proof of Corollary~\ref{cor:abstract}]
We first claim that if $K$ is a weakly compact operator in $L^1$ and $Q$ is a bounded operator in $L^1$, then $KQ$ is power compact.
Indeed we clearly have that $KQ$ is weakly compact and consequently $(KQ)^2$ is strongly compact, invoking for instance~\cite[Cor.~1 of Th.~9.9 in Chap.~II]{Schaefer-book}.
This ensures that for $S$ and $k$ as in Corollary~\ref{cor:abstract}, the operator $U=S^k$ satisfies~\eqref{eq:splittingU}.
Secondly, the existence of $\kappa_0>0$ and $g_0\in L^1_+\setminus\{0\}$ such that $Sg_0\geq \kappa_0 g_0$ readily implies, together with the positivity of $S$, that $r(U)\geq\kappa_0^k>0$.
Similarly, the existence of $\psi_0\in L^\infty_+\setminus\{0\}$ such that $S^*\psi_0\geq \kappa_0 \psi_0$ ensures that $r(U)\geq\kappa_0^k>0$.
In both cases we can thus apply Theorem~\ref{thm:U3} to $U=S^k$ in $X=L^1$.

\smallskip

We get the existence of a unique $(\phi_0,f_0)\in L^1_+\times L^\infty_+$ such that
\[Uf_0=r(U)f_0,\quad U^*\phi_0=r(U)\phi_0,\quad \text{with}\quad \|\phi_0\|_{L^\infty}=\langle\phi_0,f_0\rangle=1.\]
Applying $S$ to the first equality, we obtain $U(Sf_0)=SUf_0=r(U)Sf_0$.
Since $N(U-r(U))=\Span(f_0)$, we deduce the existence of $\Lambda_0\in\R$ such that $Sf_0=\Lambda_0f_0$.
Necessarily we have $\Lambda_0\geq0$ since $S\geq0$, and then $\Lambda_0=(r(U))^{1/k}\geq\kappa_0>0$.
We obtain in the same way that $S^*\phi_0=\Lambda_0\phi_0$.

\smallskip

Besides, since $f-\langle\phi_0,f\rangle f_0\in\{\phi_0\}^\perp$ due to the normalization $\langle\phi_0,f_0\rangle=1$, the stability result~\eqref{eq:discrete-ergo-U} implies that
\[\big\|\Lambda_0^{-kn}S^{kn}f-\langle\phi_0,f\rangle f_0\big\|_{L^1}\leq C \zeta^{kn}\left\|f-\langle\phi_0,f\rangle f_0\right\|_{L^1}\]
for all $f\in L^1$ and $n\in\N$.
Since $S$ is a bounded operator, we readily deduce the validity of~\eqref{eq:discrete-ergo} for another constant $C\geq1$.
\end{proof}

We finish by giving a Doeblin-Harris type result which is taken from~\cite[Theorem~6.3]{FGM} and allows quantifying the constants $C$ and $\zeta$ in~\eqref{eq:discrete-ergo-U}, and consequently in~\eqref{eq:discrete-ergo}, under some additional assumptions.
More precisely, we first assume that $\widetilde U$ satisfies the Lyapunov condition
\begin{equation}\label{eq:Lyapunov}
\|\widetilde Uf\|\leq\gamma\|f\|+\Theta\langle\phi_0,|f|\rangle
\end{equation}
for all $f\in X$ and some constants $\gamma\in(0,1)$ and $\Theta>0$.
Then we suppose that $\widetilde U$ satisfies the Harris condition
\begin{equation}\label{eq:Harris}
\left\{\begin{array}{l}
\exists A>\Theta/(1-\gamma),\ \exists g_A\in X_+\setminus\{0\},\ \forall f\in X_+,\vspace{2mm}\\
\|f\|\leq A\langle\phi_0,f\rangle\quad\Longrightarrow\quad \widetilde Uf\geq \langle\phi_0,f\rangle\, g_A.
\end{array}\right.
\end{equation}

\begin{thm}[\cite{FGM}]\label{thm:Harris}
Suppose that $U$ satisfies~\eqref{eq:splittingU}, \eqref{eq:Lyapunov} and \eqref{eq:Harris}.
Then the contraction estimate~\eqref{eq:discrete-ergo-U} holds true, and the constants $C$ and $\zeta$ are obtained constructively in terms of $\gamma,\Theta$ and (any positive lower bound on) the quantity $\langle\phi_0,g_A\rangle$.
\end{thm}

The proof of this result is an adapatation to the non-conservative and general Banach lattices setting of the proof of the classical Harris ergodic theorem proposed in~\cite{Canizo2023}, itself inspired from~\cite{Hairer2011} and~\cite{Gabriel2018}.

\

\section{Application to time-periodic nonlocal equations}\label{sec:EDP}

We aim at applying Corollary~\ref{cor:abstract} and then Theorem~\ref{thm:Harris} to Equation~\eqref{eq:main} for proving Theorem~\ref{thm:main}.
To this end, we need first to associate to Equation~\eqref{eq:main} a semiflow, or propagator.
In other words, we want to build a family $(S_{t,s})_{t\geq s\geq 0}$ of linear operators on $L^1(\R^d)$ such that
the unique weak solution $u(t,x)$ of~\eqref{eq:main} which start from $u_s$ at time $s$ is given by
\[u(t,x)=S_{t,s}u_s(x).\]

\medskip

The infinitesimal generator of Equation~\eqref{eq:main} is $\cL_t=\cA_t+\cB_t$ with
\[\cA_t f(x)= - \Div_x(f(x)\v(t,x)) + a(t,x)f(x)\quad\text{and}\quad \cB_t f(x)=\int_{\mathbb{R}^d} f(y) q(t,y,x)dy\]
and the backward generator of the dual equation~\eqref{eq:main-dual} is the sum of the adjoint operators
\[\cA_t^* \phi(x)= \v(t,x)\cdot\nabla_x \phi(x) + a(t,x)\phi(x)\quad\text{and}\quad \cB_t^* \phi(x)=\int_{\mathbb{R}^d} \phi(y) q(t,x,y)dy.\]
The family $(\cA_t)_{t\in\R}$ generates an explicit semiflow
\[V_{t,s}f(x)= f(X_{s,t}(x))J_{s,t}(x)e^{\int_s^t a(s',X_{s',t}(x))ds'}\]
where $X_{t,s}$ is the solution to the characteristic equation
\[\left\{\begin{array}{l}
\partial_t X_{t,s}(x)=\v(t,X_{t,s}(x)),
\vspace{2mm}\\
X_{s,s}(x)=x,
\end{array}\right.\]
and
\[J_{t,s}(x)=\det(\nabla_x X_{t,s}(x))=\exp\bigg(\int_s^t\Div_x\v(s',X_{s',s}(x))\,ds'\bigg)>0.\]
We have the useful property that for any measurable set $E\subset\R^d$,
\[|X_{s,t}(E)| = \int_{X_{s,t}(E)}dx = \int_{E}J_{s,t}(x)\,dx,\]
which ensures that
\begin{equation}\label{eq:mesureflow}
e^{-|t-s|\|\Div_x\v\|_\infty} |E| \leq |X_{s,t}(E)| \leq e^{|t-s|\|\Div_x\v\|_\infty} |E|.
\end{equation}
Denoting by $N_\infty$ the $L^\infty$ norm of $\v/(1+|x|)$, which is finite by virtue of {\bf(H$\v$)}, we have for all $x\in\R^d$ and $t\geq s$
\[|X_{t,s}(x)|+1 \leq |x| +1+ \int_{s}^{t} \displaystyle\left\lvert \v\left( u,  X_{u,s}(x) \right) \right\rvert  du \leq |x| + 1 + N_\infty\int_{s}^{t} 1+ \displaystyle\left\lvert  X_{u,s}(x) \right\rvert  du \]
and consequently, by using Gronwall's inequality,
\[|X_{t,s}(x)| +1\leq \left( |x| + 1 \right) e^{N_\infty(t-s)}.\]
Recall also that for any $t\neq s$ the mapping $x\mapsto X_{t,s}(x)$ is a $C^1$-diffeomorphism of $\R^d$ with inverse $X^{-1}_{t,s}=X_{s,t}$.
Applying the above inequality to $X_{s,t}(x)$ we can then deduce that for any $s,t\in\R$ and $x\in\R^d$
\begin{equation}\label{eq:croiss-boules}
(|x|+1)e^{-N_\infty|t-s|}-1\leq|X_{t,s}(x)|\leq (|x|+1)e^{N_\infty|t-s|}-1.
\end{equation}
The semiflow property of $(V_{t,s})$ writes
\[V_{t_3,t_2}V_{t_2,t_1}=V_{t_3,t_1}\]
for $t_3\geq t_2\geq t_1$,
and is a consequence of the fact that $u(t,x)=V_{t,s}f(x)$ is the unique weak-mild solution to the transport equation
\begin{equation*}\label{eq:transport}
\partial_t u(t,x) + \Div_x(u(t,x)\v(t,x)) = a(t,x)u(t,x)
\end{equation*}
with initial condition $u(s,x)=f(x)$.
More precisely, $(V_{t,s})$ is the unique family of operators such that
\[\langle \phi,V_{t,s}f\rangle=\langle\phi,f\rangle+\int_s^t\langle \cA^*_{s'}\phi,V_{s',s}f\rangle ds'\]
for all $t\geq s$, all $f\in L^1(\R^d)$ and all $\phi\in C^1_c(\R^d)$.
The dual semiflow writes
\[V_{s,t}^*\phi(x)= \phi(X_{t,s}(x))e^{\int_s^t a(s',X_{s',s}(x))ds'}.\]
It verifies $\langle V_{s,t}^*\phi,f\rangle = \langle \phi,V_{t,s}f\rangle$, the dual semiflow property
\[V_{s_1,s_2}^*V_{s_2,s_3}^*=V_{s_1,s_3}^*,\]
and $\varphi(s,x)=V_{s,t}^*\phi(x)$ is the solution to the dual backward equation
\begin{equation*}\label{eq:transport_dual}
-\partial_s\varphi(s,x)-\v(s,x)\cdot\nabla_x\varphi(s,x)= a(s,x)\varphi(s,x).
\end{equation*}
with terminal condition $\varphi(t,x)=\phi(x)$.

\

The operator $\cB_t$ is a bounded operator, and $\cL_t$ is thus a bounded perturbation of $\cA_t$ which generates a positive semiflow $(S_{t,s})$ that verifies the Duhamel formulae
\begin{equation}\label{eq:Duhamel1}
S_{t,s} = V_{t,s} + \int_s^t S_{t,s'}\cB_{s'}V_{s',s}\,ds',
\end{equation}
\begin{equation}\label{eq:Duhamel2}
S_{t,s} = V_{t,s} + \int_s^t V_{t,s'}\cB_{s'}S_{s',s}\,ds',
\end{equation}
and their dual counterparts
\begin{equation}\label{eq:Duhamel1_dual}
S_{s,t}^* = V_{s,t}^* + \int_s^t V_{s,s'}^*\cB^*_{s'}S_{s',t}^*\,ds',
\end{equation}
\begin{equation}\label{eq:Duhamel2_dual}
S_{s,t}^* = V_{s,t}^* + \int_s^t S_{s,s'}^*\cB^*_{s'}V_{s',t}^*\,ds'.
\end{equation}
From these formula we easily infer, by using Grönwall's inequality and the Assumptions~{\bf(H$a$)} and~\eqref{as:qsup}, the growth estimate
\begin{equation}\label{eq:expo-growth}
\|S_{t,s}\|=\|S^*_{s,t}\|\leq e^{(\hat q+\esssup a)(t-s)}.
\end{equation}
It is standard result that $(V_{t,s})$ is a strongly continuous semiflow in $L^1(\R^d)$, and consequently so does $(S_{t,s})$.
The semiflow $(S_{t,s})$ yields the unique weak-mild solutions of Equation~\eqref{eq:main}, in the sense that
\begin{equation}\label{eq:weak-mild}
\langle \phi,S_{t,s}f\rangle=\langle\phi, f\rangle+\int_s^t\langle \cL^*_{s'}\phi,S_{s',s}f\rangle ds'
\end{equation}
for all $t\geq s\geq0$, all $f\in L^1(\R^d)$ and all $\phi\in C^1_c(\R^d)$.
The dual semiflow $(S^*_{s,t})$ yields the solutions to Equation~\eqref{eq:main-dual}, and the dual counterpart of~\eqref{eq:weak-mild} reads
\begin{equation}\label{eq:weak-mild-dual}
\langle S^*_{s,t}\phi,f\rangle=\langle\phi,f\rangle+\int_s^t\langle S^*_{s',t}\phi,\cL_{s'}f\rangle ds'
\end{equation}
for all $t\geq s\geq0$, all $\phi\in L^\infty(\R^d)$ and all $f\in C^1_c(\R^d)$.
The $T$-periodicity of $\cL_t$ leads to the following $T$-periodicity property of the semiflow
\[S_{t+T,s+T}=S_{t,s}.\]

\medskip

We now define the bounded operator $S=S_{T,0}$, and for proving Theorem~\ref{thm:main}, we verify that it satisfies the conditions in Corollary~\ref{cor:abstract}, as well as the conditions allowing to apply Theorem~\ref{thm:Harris}.
To do so, we start with some preliminary results.

\subsection{Preliminary lemmas}

\begin{lem}\label{lem:g_0}
Under Assumptions~{\bf(H$\v$)}-{\bf(H$a$)}-{\bf(H$q$)}, there exists $\kappa_0>0$ and $g_0\in L^1_+\setminus\{0\}$ such that $Sg_0\geq\kappa_0 g_0$.
\end{lem}

\begin{proof}
We use the positivity condition~\eqref{as:qpos2} and we consider a function $g_0 \in C^1( \mathbb{R}^d)$ such that $g_0>0$ on $B(x_0,r_0)$ and $g_0=0$ on $\R^d  \setminus B(x_0,r_0)$.
For any $r\in(r_1,r_0)$ fixed (to be chosen later), $x \in B(x_0,r)$, and $t \geq 0$, we have
\begin{align*}
   \cL_t g_0(x) &=  - \Div_x(g_0(x)\v(t,x)) + a(t,x)g_0(x)+\int_{\R^d} g_0(y) q(t,y,x)dy \\
   &\geq - \left(\max_{\overline B(x_0,r)}|\nabla\!_x g_0| \right) \left( \esssup_{(0,T)\times B(x_0,r)} |\v|\right) - \left(\esssup_{(0,T)\times B(x_0,r)} |a -\Div_x\v|\right) g_0(x)\\
   &\geq -\left(\frac{ \max_{\overline B(x_0,r)} |\nabla\!_x g_0| }{ \min_{\overline B(x_0,r)} |g_0|}  \esssup_{(0,T)\times\overline B(x_0,r)}|\v|  + \esssup_{(0,T)\times B(x_0,r)} |a -\Div_x\v| \right) g_0(x).
\end{align*}

\medskip

\noindent For $x \in A_r:=B(x_0,r_0) \setminus B(x_0,r)$ we have
\begin{align*}
       \cL_t g_0(x)& =  - \Div_x(g_0(x)\v(t,x)) + a(t,x)g_0(x)+\int_{\R^d} g_0(y) q(t,y,x)dy \\
       &\geq - \left(\sup_{A_r} |\nabla\!_x g_0 |  \right)\! \left( \esssup_{(0,T)\times B(x_0,r_0)} \lvert \v \rvert \right) - \left(\esssup_{(0,T)\times B(x_0,r_0)} |a -\Div_x\v|\right) g_0(x) + q_0\inf_{B(x_0,r_1)}g_0.
\end{align*}
Since $\sup_{A_r}|\nabla\!_x g_0 |\to0$ when $r\to r_0$, we can find $r$ close enough to $r_0$ so that
\[- \left(\sup_{A_r} |\nabla\!_x g_0 |  \right)\! \left( \esssup_{(0,T)\times B(x_0,r_0)} \lvert \v \rvert \right) + q_0\inf_{B(x_0,r_1)}g_0 \geq 0\]
and we thus have $\cL_t g_0(x)\geq-\alpha_0g_0(x)$ for all $x \in \R^d$ and all $t\in\R$ with
\[\alpha_0 = \frac{ \max_{\overline B(x_0,r)} |\nabla\!_x g_0| }{ \min_{\overline B(x_0,r)} |g_0|}  \esssup_{(0,T)\times\overline B(x_0,r)}|\v|  + \esssup_{(0,T)\times B(x_0,r_0)} |a -\Div_x\v|.\]
Injecting this inequality in~\eqref{eq:weak-mild-dual} with $f=g_0$ we obtain that for any $\phi\geq0$
\[\langle \phi,S_{t,s}g_0\rangle=\langle S^*_{s,t}\phi,g_0\rangle\geq\langle\phi,g_0\rangle-\alpha_0\int_s^t\langle S^*_{s',t}\phi,g_0\rangle ds'\geq\langle\phi,g_0\rangle-\alpha_0\int_s^t\langle \phi,S_{t,s'}g_0\rangle ds'.\]
The Grönwall inequality applied to the function $s\mapsto\langle \phi,S_{T,s}g_0\rangle$ then gives that
\[\langle \phi,S_{T,s}g_0\rangle\geq\langle \phi,g_0\rangle e^{-\alpha_0(T-s)}\]
for all $\phi\geq0$ and all $s\in[0,T]$.
Taking $s=0$ we deduce that $Sg_0\geq e^{-\alpha_0T}g_0$ and the proof is complete.
\end{proof}

\medskip

\begin{lem}\label{lem:U=W+K}
Assume that~{\bf(H$\v$)}-{\bf(H$a$)}-{\bf(H$q$)} are satisfied.
Then for any $\kappa>0$, there exists an integer $k$ such that $S^k$ can be split as $S^k=W+K$ with $W$ such that $\|W\|<\kappa^k$ and $K$ weakly compact.
\end{lem}

\begin{proof}
We start from the Duhamel formula~\eqref{eq:Duhamel2} to write, for any $k\in\N$ and $\delta\in(0, kT)$
\[S^k = V_{kT,0} + \int_{kT-\delta}^{kT}V_{kT,s}\cB_{s}S_{s,0}\,ds + \int_0^{kT-\delta} V_{kT,s}\cB_{s}S_{s,0}\,ds.\]
In a first step we prove by using the Dunford-Pettis theorem that the last integral term is a weakly compact operator in $L^1$ for any $k\in\N$ and $\delta\in(0, kT)$.
In a second step, we prove that we can find $k$ large enough to have $\|V_{kT,0}\|<\kappa^k$, which then allows choosing $\delta$ small enough so that $\big\|V_{kT,0} + \int_{kT-\delta}^{kT}V_{kT,s}\cB_{s}S_{s,0}\,ds\big\|<\kappa^k$.

\medskip

\paragraph{\it Step~1.}
Let us fix any integer $k$ and any $\delta\in(0,kT)$, and set
\[K=\int_0^{kT-\delta} V_{kT,s}\cB_{s}S_{s,0}\,ds.\]
We want to prove that $K$ is a weakly compact operator by invoking the Dunford-Pettis theorem.
To do so we consider a sequence $(f_n)_{n\geq0}$ such that $f_n\geq0$ and $\|f_n\|_{L^1}\leq1$ for any $n\geq0$, and we take $\epsilon>0$.
We need to prove the existence of $\eta>0$ and $R>0$ such that for all set $E\subset\R^d$ we have
\[\left.\begin{array}{r}
|E|<\eta\\
\text{or}\\
E\subset B^c(0,R)
\end{array}\right\}
\quad\Longrightarrow\quad \int_EKf_n<\epsilon.\]
Let $E$ be a measurable subset of $\R^d$.
We observe that for any $n\geq0$ and any $s\in[0,kT-\delta]$
\[\int_{E} V_{kT,s}\cB_{s}S_{s,0} f_n(x) \,dx = \langle \mathbf{1}_E , V_{kT,s} \cB_{s} S_{s,0}f_n \rangle = \langle \cB^*_{s} V^*_{s,kT}  \mathbf{1}_E , S_{s,0} f_n \rangle \leq   {\| \cB^*_{s} V^*_{s,kT}  \mathbf{1}_E   \|}_{L^\infty} {\| S_{s,0} f_n \|}_{L^1}.\]
Setting $\omega:=\max(1,\hat q+\esssup\overline a)>0$ we deduce from~\eqref{eq:expo-growth} that for all $n\geq0$
\[\int_EKf_n\leq \frac{1}{\omega}e^{\omega kT}\sup_{0\leq s\leq kT-\delta}{\| \cB^*_{s} V^*_{s,kT}  \mathbf{1}_E   \|}_{L^\infty},\]
and we have for all $x\in\R^d$
\[\cB^*_{s} V^*_{s,kT}  \mathbf{1}_E(x)  = \int_{\mathbb{R}^d} q(s,x,y) \mathbf{1}_E \left( X_{kT,s}(y) \right) e^{ \int_{s}^{kT} a(s', X_{s',s}(y)) ds' } dy \leq e^{kT \esssup \overline{a}} \int_{X_{kT,s}^{-1}(E)}q(s,x,y)\,dy.\]
From \eqref{cond Q} we get the existence of $\eta_0>0$ such that for any $E_0\subset\R^d$
\[|E_0| < \eta_0 \quad \Longrightarrow\quad \sup_{t}\esssup_{x} \int_{E_0} q(t,x,y)dy < \epsilon\,\omega e^{-kT (\omega+\esssup \overline{a})}.\]
As recalled in~\eqref{eq:mesureflow}, we have
\[|X_{kT,s}^{-1}(E)|=|X_{s,kT}(E)| \leq e^{(kT-s)\|\Div_x\v\|_\infty} |E|\]
so by choosing $\eta =\eta_0e^{-kT\|\Div_x\v\|_\infty}$ we obtain that if $|E|<\eta$ then for all $n\geq0$
\[\int_EKf_n  < \epsilon.\]
Now we consider $E \subset  \mathbb{R}^d \setminus B(0,R)$.
Due to the assumption {\bf(H$a$)}, we can define the radial function $\overline a$ by
\[\overline a(x)=\esssup_{t\in\R,|y|>|x|} a(t,y)=\esssup_{t\in[0,T],|y|>|x|} a(t,y),\]
which verifies
\[\lim_{|x|\to\infty}\overline a(x)=-\infty\qquad \text{and}\qquad a(t,x)\leq\overline a(x)\quad\text{for}\ a.e.\ (t,x)\in\R\times\R^d.\]
In particular, there exists $R_a>0$ such that
\[\esssup_{B^c(0,R_a)}\overline a<0.\]
We thus have for any $x\in\R^d$ and $s\in[0,kT-\delta]$
\begin{align*}
    \cB^*_{s} V^*_{s,kT}  \mathbf{1}_E(x)  &= \int_{\R^d} q(s,x,y) \mathbf{1}_E ( X_{kT,s}(y)) e^{ \int_{s}^{kT} a(s', X_{s',s}(y)) ds' }dy\\
    &  \leq \int_{X_{s,kT}(E)} q(s,x,y)e^{ \int_{s}^{kT} \overline{a}( X_{s',s}(y)) ds' }dy\\
    & \leq \bigg(\sup_{s\in[0, kT-\delta]}\esssup_{y\in X_{s,kT}(E)}e^{ \int_{s}^{kT} \overline{a}( X_{s',s}(y)) ds'}\bigg) \hat q. 
\end{align*}
For $y\in X_{s,kT}(E)$, there is $x\in E\subset B^c(0,R)$ such that $y=X_{s,kT}(x)$, and so $X_{s',s}(y)= X_{s',kT}(x)$.
We deduce from~\eqref{eq:croiss-boules} that
\[|X_{s',s}(y)|=|X_{s',kT}(x)|\geq |x|e^{-N_\infty(kT-s)}-1\geq Re^{N_\infty kT}-1.\]
Choosing $R$ large enough so that $Re^{N_\infty kT}-1\geq R_a$ and accordingly $\esssup_{B^c(0,Re^{N_\infty kT}-1)}\overline a<0$, we obtain that
\begin{align*}
\sup_{s\in[0, kT-\delta]}\esssup_{y\in X_{s,kT}(E)}e^{ \int_{s}^{kT} \overline{a}( X_{s',s}(y)) ds'}&
\leq \sup_{s\in[0, kT-\delta]}\exp\Big((kT-s)\,\esssup_{B^c(0,Re^{N_\infty kT}-1)}\overline a\Big)\\
&\leq \exp\Big(\delta\,\esssup_{B^c(0,Re^{N_\infty kT}-1)}\overline a\Big).
\end{align*}
The right hand side tends to zero when $R\to+\infty$, because of~\eqref{as:asup}, so we can find $R$ large enough so that $\sup_{0\leq s\leq kT+\delta}{\| \cB^*_{s} V^*_{s,kT}  \mathbf{1}_E   \|}_{L^\infty}<\epsilon\,\omega e^{-\omega kT}$ and consequently, for all $n\geq0$,
\[\int_EKf_n  < \epsilon.\]

\medskip

\paragraph{\it Step~2.}
We want to find $k$ large enough and $\delta$ small enough so that the norm of $W$, defined by
\[W = V_{kT,0} + \int_{kT-\delta}^{kT}V_{kT,s}\cB_{s}S_{s,0}\,ds,\]
satisfies $\|W\|<\kappa^k$.
We observe that for all $f\in L^1$
\[\| V_{kT,0}f  \|_{L^1} =  \sup_{ \| \phi \|_{L^\infty} \leq 1}\langle f , V^*_{0,kT}\phi \rangle
\leq \int_{\R^d}|f(x)|  e^{\int_{0}^{kT} \overline{a}(X_{s,0}(x))ds } dx,\]
so that
\[\| V_{kT,0}\| \leq \sup_{x\in\R^d} e^{\int_{0}^{kT} \overline{a}(X_{s,0}(x))ds }.\]
Defining for any $R>0$ and any $x\in\R^d$ the set
\[\Delta_R(x) := \big\{ t  \geq 0, \ X_{t,0}(x) \in B(0,R) \big\},\]
the condition~\eqref{eq:DeltaR} reads
\[ \forall R>0, \quad \underset{x \in \mathbb{R}^d}{\sup} \lvert \Delta_R(x) \rvert  <\infty.\]
We can thus write, for any $k \in\N$, $x\in\R^d$, and $R>R_a$, so that $\esssup_{B^c(0,R)}\overline{a}<0$,
\begin{align*}
    \int_{0}^{kT} \overline{a}(X_{s,0}(x))ds &= \int_{[0,kT]\cap \Delta_R(x)}\overline{a}(X_{s,0}(x))ds + \int_{[0,kT] \cap\Delta^c_R(x)}^{} \overline{a}(X_{s,0}(x))ds \\
    &\leq  \sup_{x \in \R^d}|\Delta_R(x)| \bigg(  \esssup_{B(0,R)}\overline{a}  \bigg) + \bigg(  \esssup_{B^c(0,R)}\overline{a}  \bigg)\Big(kT-\underset{x \in \R^d}{\sup}|\Delta_R(x)|\Big) \\
    &\leq  \bigg(  \esssup_{B^c(0,R)}\overline{a}  \bigg)kT + \sup_{x \in \R^d}|\Delta_R(x)| \bigg(  \esssup_{B(0,R)}\overline{a} - \esssup_{B^c(0,R)}\overline{a}  \bigg).
\end{align*}
We suppose without loss of generality that $\kappa\in(0,1)$ and we choose $R$ large enough so that
\[\esssup_{B^c(0,R)}\overline{a} \leq \frac{\log\kappa}{T}-1.\]
Consequently
\[\sup_{x\in\R^d}\int_{0}^{kT} \overline{a}(X_{s,0}(x))ds\leq k\log\kappa-kT+ \sup_{x \in \R^d}|\Delta_R(x)| \bigg(  \esssup_{B(0,R)}\overline{a} - \esssup_{B^c(0,R)}\overline{a}  \bigg)\]
and we can take $k$ large enough to have
\[\sup_{x\in\R^d}\int_{0}^{kT} \overline{a}(X_{s,0}(x))ds\leq k\log\kappa-1,\]
which yields $\| V_{kT,0}\| \leq \kappa^k e^{-1}$.
Finally, using that the function $s\mapsto\|V_{kT,s}\cB_{s}S_{s,0}\|$ is bounded on $[0,kT]$, we can find $\delta>0$ small enough so that
\[\bigg\|\int_{kT-\delta}^{kT}V_{kT,s}\cB_{s}S_{s,0}\,ds\bigg\| < \kappa^k(1-e^{-1}),\]
and we thus obtain $\|W\|<\kappa^k$.
\end{proof}

\medskip

\begin{lem}\label{lem:positivity}
Assume that~{\bf(H$\v$)}-{\bf(H$a$)}-{\bf(H$q$)} are satisfied.
Then for all $t>s\geq0$, the operator $S_{t,s}$ is strongly positive, in the sense that
\[\forall f\in L^1_+\setminus\{0\},\quad S_{t,s}f(x)>0\quad \text{for almost every}\ x\in\R^d.\]
\end{lem}

\begin{proof}
Let $f\geq0$ such that $f\neq0$ and let $t>s\geq0$.
Since $S_{t,s}f\geq V_{t,s}f$ we also have that $S_{t,s}f\geq0$ and $S_{t,s}f\neq0$, so that $\suppess S_{t,s}f\neq\emptyset$. We consider now $x\in\suppess S_{t,s}f$ and we prove that $S_{t,s}f>0$ almost everywhere on $B(x,r_0/2)$, where $r_0>0$ is defined in Assumption~\eqref{as:qpos}.
Let $E\subset B(x,r_0/2)$ such that $S_{t,s}f=0$ almost everywhere on $E$.
By positivity, we deduce from the Duhamel formula~\eqref{eq:Duhamel2} that
\[\int_s^t \bigg(\int_E V_{t,s'}\cB_{s'}S_{s',s}f\bigg)ds'=0,\]
by using the Fubini-Tonelli theorem.
We easily check that $s'\mapsto\cB_{s'}$ is a strongly continuous family of operators in $L^1$, by using Assumption~\eqref{as:qsup} and the time continuity of $q$.
We deduce that the function $s'\mapsto \int_E V_{t,s'}\cB_{s'}S_{s',s}f$ is continuous.
Since it is non-negative and its integral on $[s,t]$ is null, it must be identically zero, and in particular we get by considering $s'=t$ that
\[\int_E\cB_tS_{t,s}f=0,\]
and a fortiori
\[\int_E\int_{B(x,r_0/2)}q(t,y,z)(S_{t,s}f)(y)\,dydz=0.\]
Since $|z-y|<r_0$ for any $z\in E\subset B(x,r_0/2)$ and $y\in B(x,r_0/2)$, we infer from positivity assumption~\eqref{as:qpos} that
\[\int_E \int_{B(x,r_0/2)}(S_{t,s}f)(y)\,dydz=0.\]
This enforces $E$ to be a set of measure zero since $x\in\suppess S_{t,s}f$ implies that $\int_{B(x,r_0/2)}S_{t,s}f>0$.

\smallskip

We have proved that $S_{t,s}f>0$ almost everywhere in $B(x,r_0/2)$ for any $x\in\suppess S_{t,s}f$. This ensures in particular that $\suppess S_{t,s}f$ is open, and so $\suppess S_{t,s}f=\R^d$.
We have thus proved that $S_{t,s}f>0$ almost everywhere in $B(x,r_0/2)$ for any $x\in\R^d$, or in other words, that $S_{t,s}f>0$ almost everywhere in~$\R^d$.
\end{proof}

\medskip

Lemmas~\ref{lem:g_0}, \ref{lem:U=W+K}, and~\ref{lem:positivity} allow applying Corollary~\ref{cor:abstract} to the operator $S$ and thus proving the non-constructive part of Theorem~\ref{thm:main}.
For quantifying the constants $C$ and $\rho$ by invoking Theorem~\ref{thm:Harris}, we need some additional estimates that we will prove now.

\medskip

\begin{lem}\label{lem:Harris0}
Assume that~{\bf(H$\v$)}-{\bf(H$a$)}-{\bf(H$q+$)} are satisfied.
Then for any $t>s\geq0$ and any $r,R>0$, there exists $c>0$ such that
\[S_{s,t}^*\1_{B(0,r)}\geq c\,\1_{B(0,R)}\qquad\text{and}\qquad S_{t,s}\1_{B(0,r)}\geq c\,\1_{B(0,R)}.\]
\end{lem}

\begin{proof}
Let $r>0$.
Due to~\eqref{eq:croiss-boules}, we can find $\epsilon_r>0$ such that
\begin{equation}\label{eq:croiss-boules-2}
|t-s|<\epsilon_r\quad\Longrightarrow\quad|X_{t,s}(x)|<|x|+\frac{r_0}{8},\quad\forall x\in B\Big(0,r+\frac{r_0}{2}\Big),
\end{equation}
where $r_0$ is defined in~{\bf(H$q+$)}.
Note that, up to changing $r_0$ into $\min(r,r_0)$ in~{\bf(H$q+$)}, we can assume without loss of generality that $r_0\leq r$.
Denoting $a_r:=\essinf_{[0,T]\times B(0,r+r_0)}a$ and iterating once the Duhamel formula~\eqref{eq:Duhamel1_dual} or~\eqref{eq:Duhamel2_dual}, we can write by using {\bf(H$q+$)} that for any $t\geq s\geq0$ such that $|t-s|<\epsilon_r$ and any $x\in B(0,r+r_0/2)$
\begin{align*}
S^*_{s,t}\1_{B(0,r)}(x)&\geq \int_s^t V^*_{s,s'}\cB^*_{s'}V^*_{s',t}\1_{B(0,r)}(x)\,ds'\\
&\geq \int_s^t e^{(t-s')a_r}V^*_{s,s'}\cB^*_{s'}\1_{B(0,r)}(X_{t,s'}(x))\,ds'\\
&\geq \int_s^t e^{(t-s')a_r}V^*_{s,s'}\int_{\R^d}q(s',x,y)\1_{B(0,r)}(X_{t,s'}(y))\,dy\,ds'\\
&\geq q_0\int_s^t e^{(t-s')a_r}V^*_{s,s'}\int_{B(x,r_0)}\1_{B(0,r)}(X_{t,s'}(y))\,dy\,ds'\\
&\geq e^{(t-s)a_r}q_0\int_s^t\int_{B(X_{s',s}(x),r_0)}\1_{B(0,r)}(X_{t,s'}(y))\,dy\,ds'.
\end{align*}
Invoking~\eqref{eq:croiss-boules-2} we obtain
\[S^*_{s,t}\1_{B(0,r)}(x)\geq e^{(t-s)a_r}q_0\int_s^t\int_{B(X_{s',s}(x),r_0)}\1_{B(0,r-r_0/8)}(y)\,dy\,ds'.\]
Using again~\eqref{eq:croiss-boules-2}, we have that $X_{s',s}(x)\in B(0,r+\frac{5r_0}{8})$ for any $s'\in[s,t]$ and any $x\in B(0,r+r_0/2)$.
It is clear due to radial symmetry that the measure of the intersection $B\big(0,r-\frac{r_0}{8}\big)\cap B\big(x_r,r_0\big)$ is the same for any $x_r\in\R^d$ such that $|x_r|=r+\frac{5r_0}{8}$, and that it is positive.
We thus set
\[b_r:=\Big|B\Big(0,r-\frac{r_0}{8}\Big)\cap B\big(x_r,r_0\big)\Big|>0\]
for an arbitrary choice of $x_r$ satisfying $|x_r|=r+\frac{5r_0}{8}$ and we claim that for all $s'\in[s,t]$ and all $x\in B(0,r+r_0/2)$
\[\Big| B\Big(0,r-\frac{r_0}{8}\Big)\cap B\big(X_{s',s}(x),r_0\big)\Big|\geq b_r.\]
To verify it, we define
\[\alpha(s',x)=\frac{|X_{s',s}(x)|}{r+\frac{5r_0}{8}}\in[0,1)\qquad\text{and}\qquad x_r(s',x)=
\begin{cases}\frac{1}{\alpha(s',x)}X_{s',s}(x)&\ \text{if}\ \alpha(s',x)>0\\
x_r&\ \text{if}\ \alpha(s',x)=0\end{cases},\]
so that $|x_r(s',x)|=r+\frac{5r_0}{8}$ and $X_{s',s}(x)=\alpha(s',x) x_r(s',x)$.
Since it is clear that the mapping $m:[0,1]\to\R_+$ defined by
\[m(\alpha)=\Big| B\Big(0,r-\frac{r_0}{8}\Big)\cap B\big(\alpha x_r(s',x),r_0\big)\Big|\]
is nonincreasing, we have that $m(\alpha(s',x))\geq m(1)=b_r$, which is the claim.
We have thus proved that for any $t\geq s\geq0$ such that $|t-s|<\epsilon_r$
\[S^*_{s,t}\1_{B(0,r)}\geq q_0b_r(t-s)e^{(t-s)a_r}\1_{B(0,r+r_0/2)}.\]
Setting $r'=r+\frac{r_0}{2}$ and taking $t'\in[s,t]$ such that $|t'-s|<\epsilon_{r'}$ we get
\begin{align*}
S^*_{s,t}\1_{B(0,r)}=S^*_{s,t'}S_{t',t}\1_{B(0,r)}&\geq q_0b_r(t-t')e^{(t-t')a_r}S_{s,t'}\1_{B(0,r+r_0/2)}\\
&\geq q_0^2b_rb_{r'}(t-t')(t'-s)e^{(t-t')a_r+(t'-s)a_{r'}}\1_{B(0,r+r_0)}.
\end{align*}
Iterating this argument, we obtain that for any $r,R>0$ and any $t>s\geq0$ such that $|t-s|<\epsilon_r$ we have
\[S^*_{s,t}\1_{B(0,r)}\geq c\,\1_{B(0,R)}\]
for some constant $c=c(r,R,s,t)>0$.
Splitting the time interval $[s,t]$ into small enough sub-intervals $[t_i,t_{i+1}]$, we readily infer that this inequality still holds without the restriction $|t-s|<\epsilon_r$.

\medskip

The proof for $S_{t,s}$ follows exactly the same arguments and we then skip it.
\end{proof}

\medskip

\begin{cor}\label{cor:Harris}
Assume that~{\bf(H$\v$)}-{\bf(H$a$)}-{\bf(H$q+$)} are satisfied.
Then for any $R_1,R_2>0$, there exists $c>0$ such that for all $\phi\in L^\infty_+$
\begin{equation}\label{eq:HarrisR1R2}
S^*\phi\geq c\,\langle\phi,\1_{B(0,R_1)}\rangle\1_{B(0,R_2)}.
\end{equation}
\end{cor}

\begin{proof}
Similarly as in the proof of Lemma~\ref{lem:Harris0}, we can find $\epsilon>0$ small enough so that
\begin{equation}\label{eq:croiss-boules-3}
|t-s|<\epsilon\quad\Longrightarrow\quad|X_{t,s}(x)|<|x|+\frac{r_0}{4},\quad\forall x\in B(0,r_0),
\end{equation}
where $r_0$ is still defined in {\bf(H$q+$)},
and then the existence of $c_\epsilon>0$ such that for any $\phi\in L^\infty_+$ and any $x\in B(0,r_0/4)$
\[S^*_{s,s+\epsilon}\phi(x)\geq c_\epsilon\int_s^{s+\epsilon}\int_{B(X_{s',s}(x),r_0)}\phi(X_{s+\epsilon,s'}(y))\,dy\,ds'.\]
Using~\eqref{eq:croiss-boules-3} and the fact that $x\in B(0,r_0/4)$ we infer that
\[S^*_{s,s+\epsilon}\phi(x)\geq c_\epsilon\int_s^{s+\epsilon}\int_{B(0,r_0/2)}\phi(X_{s+\epsilon,s'}(y))\,dy\,ds'.\]
Through a change of variable we obtain
\[S^*_{s,s+\epsilon}\phi(x)\geq c_\epsilon \int_s^{s+\epsilon}\!\!e^{-(s+\epsilon-s')\|\Div_x\v\|_\infty}\!\! \int_{X_{s',s+\epsilon}(B(0,r_0/2))}\phi(z)\,dzds'\geq \epsilon c_\epsilon e^{-\epsilon\|\Div_x\v\|_\infty}\!\! \int_{B(0,r_0/4)}\phi(z)\,dz,\]
since~\eqref{eq:croiss-boules-3} guarantees that $X_{s+\epsilon,s'}(B(0,r_0/4))\subset B(0,r_0/2)$, which also reads $B(0,r_0/4)\subset X_{s',s+\epsilon}(B(0,r_0/2))$.
Choosing $s=T/2$, we have proved the existence of $\epsilon\in(0,T/2)$ and $\tilde c_\epsilon>0$ such that for any $\phi\geq0$
\[S^*_{\frac{T}{2},\frac{T}{2}+\epsilon}\phi\geq \tilde c_\epsilon \langle\phi,\1_{B(0,r_0/4)}\rangle\1_{B(0,r_0/4)}.\]
Using Lemma~\ref{lem:Harris0} we finally get that for any $R_1,R_2>0$ and any $\phi\geq0$
\[
S^*\phi=S^*_{0,\frac T2}S^*_{\frac T2,\frac T2+\epsilon}S^*_{\frac T2+\epsilon,T}\phi
\geq \tilde c_\epsilon \langle \phi,S_{T,\frac T2+\epsilon}\1_{B(0,r_0/4)}\rangle S^*_{0,\frac T2}\1_{B(0,r_0/4)}
\geq c\,\langle\phi,\1_{B(0,R_1)}\rangle\1_{B(0,R_2)}
\]
for some $c>0$.
\end{proof}

\medskip

\begin{lem}\label{lem:LinfL1}
Assume that~{\bf(H$\v$)}-{\bf(H$a$)}-{\bf(H$q$)} are satisfied.
Then for any $\lambda>\|W\|$, 
there exists $\alpha<1$ and $h \in L^1_+\setminus\{0\}$ such that for all $\phi \in L^\infty$
\begin{equation}\label{eq:LinfL1}
{\| (W^*-\lambda)^{-1}K^* \phi \|}_{L^\infty} \leq \alpha {\| \phi \|}_{L^\infty}+\langle \phi, h \rangle.
\end{equation}
\end{lem}

\begin{proof}
For $\phi \in L^\infty$, we have for all $\phi\in L^\infty$ and all $x\in\R^d$
\[ K^*\phi(x)  =  \int_{0}^{kT-\delta} S^*_{0,s}  \cB^*_{s} V_{s,kT}^* \phi(x) ds 
=  \int_{0}^{kT-\delta} S^*_{0,s}\left( \int_{\R^d} e^{\int_{s}^{kT} a(s',X_{s',s}(y))ds'} \phi\left(X_{kT,s}(y) \right) q(s,y,x)dy \right)ds.\]
Setting $\omega:=\max(1,\hat q+\esssup\overline a)>0$ we deduce from~\eqref{eq:expo-growth} that
\[ {\|K^*\phi\|}_{L^\infty} \leq \frac{e^{\omega kT}}{\omega}\sup_{s\in[0,kT-\delta]}\esssup_{x\in\R^d}\left( \int_{\R^d} e^{\int_{s}^{kT} a(s',X_{s',s}(y))ds'} |\phi(X_{kT,s}(y))| \,q(s,y,x)dy \right).\]
We split the integral on $\R^d$ as
\begin{align*}
 \int_{\mathbb{R}^d} e^{\int_{s}^{kT} a(s',X_{s',s}(y))ds'}& \phi\left(X_{kT,s}(y) \right) q(s,y,x)dy  \\
 &= \int_{B(0,R)} e^{\int_{s}^{kT} a(s',X_{s',s}(y))ds'} q(s,y,x)  \phi\left(X_{kT,s}(y) \right) dy ds\\
 &\hskip20mm + \int_{B^c(0,R)} e^{\int_{s}^{kT} a(s',X_{s',s}(y))ds'} q(s,y,x)  \phi\left(X_{kT,s}(y) \right) dy ds.
\end{align*}
Similarly as in the proof of Lemma~\ref{lem:U=W+K}, we can choose $R$ large enough in order to have
$$\int_{B^c(0,R)} e^{\int_{s}^{kT} a(s',X_{s',s}(y))ds'} q(s,y,x) \phi\left(X_{kT,s}(y) \right) dy ds \leq \frac {\omega e^{-\omega kT}\|\phi\|_{L^\infty}}{4 \| (W^*-\lambda)^{-1} \|}$$
for all $s\in[0,kT-\delta]$ and all $x\in\R^d$.
Let us now fix a representative of $q$ in the quotient space $C(\R,L^\infty(\R^d,L^1(\R^d,\R))$
and define for any $s\in\R$, $x\in\R^d$ and $M>0$ the set
\[E_M(s,x)=\{ y\in\R^d, q(s,x,y)\geq M\}.\]
We then have
\begin{align*}
    \int_{B(0,R)} & e^{\int_{s}^{kT} a(s',X_{s',s}(y))ds'} q(s,y,x) \phi\left(X_{kT,s}(y) \right) dy \\
    &= \int_{B(0,R)\cap E_M(s,x)} e^{\int_{s}^{kT} a(s',X_{s',s}(y))ds'} q(s,y,x) \phi\left(X_{kT,s}(y) \right) dy\\
    &\hspace{24mm} + \int_{B(0,R)\cap E^c_M(t,x)} e^{\int_{s}^{kT} a(s',X_{s',s}(y))ds'} q(s,y,x) \phi\left(X_{kT,s}(y) \right) dy\\
    &\leq  e^{kT \esssup \overline{a}} \left\lVert \phi \right\lVert_{L^\infty}  \int_{B(0,R)\cap E_M(s,x)} q(s,y,x) dy + M e^{kT \esssup \overline{a}} \int_{\R^d}\mathbf{1}_{B(0,R)}(y) \phi\left(X_{kT,s}(y) \right) dy \\
    &= e^{kT \esssup \overline{a}} {\|\phi\|}_{L^\infty}\int_{B(0,R)\cap E_M(s,x)} q(s,y,x) dy +   M e^{kT \esssup \overline{a}}  \int_{\R^d} J_{s,kT}(z) \mathbf{1}_{B(0,R)}\left( X_{s,kT}(z) \right)  \phi(z)  dz.
\end{align*}
From~\eqref{eq:croiss-boules} we have $|X_{t,s}(z)| \leq \left( |z| + 1 \right) e^{ N_\infty(t-s)} $ and consequently
$$ \mathbf{1}_{B(0,R)}  \left( X_{s,kT}(z) \right) \leq  \mathbf{1}_{B\left( 0,(R+1)e^{N_\infty kT} \right)}(z) .$$
From (\ref{cond Q}) there exist $\eta >0$ such that $$\quad \forall \lvert E \rvert < \eta, \quad \sup_t\,\esssup_y \displaystyle \int_{E} q(t,x,y)dy < \frac{\omega e^{-(\esssup\overline a+\omega)kT}}{4 \| (W^*-\lambda)^{-1} \|}.$$
We have from Markov's inequality $$ |E_M(s,x)| \leq \frac{1}{M} \int_{\mathbb{R}^d}^{} q(s,x,y)dy .$$
We take $M$ large enough in order to have $  \underset{s \in[0,kT]}{\sup}\underset{x }{\esssup} \lvert B(0,R)\cap E_M(s,x) \rvert< \eta $ and we finally obtain 
$$  \left\lVert(W^*-\lambda)^{-1} K^*\phi \right\rVert_{L^\infty}\leq\|(W^*-\lambda)^{-1}\|\,{\|K^*\phi\|}_{L^\infty} \leq \frac{1}{2}\left\lVert \phi \right\lVert_{L^\infty} + \langle \phi, h \rangle    $$
where $h$ is defined by
\begin{equation}\label{eq:h}
h(z) = \|(W^*-\lambda)^{-1}\| M e^{kT \esssup \overline{a}}  \left( \underset{s \in[0,kT]}{\sup}  J_{s,kT}(z) \right) \mathbf{1}_{B\left( 0,(R+1)e^{N_\infty kT} \right)}(z).
\end{equation}
Clearly, changing the representative of $q$ does not affect this inequality and the lemma is proved.
\end{proof}

\medskip

We have now all the ingredients for proving Theorem~\ref{thm:main}

\subsection{Proof of Theorem~\ref{thm:main}}

We start by verifying the conditions in Corollary~\ref{cor:abstract} under Assumptions~{\bf(H$\v$)}-{\bf(H$a$)}-{\bf(H$q$)}. 
The condition~1 is exactly Lemma~\ref{lem:g_0}, applying Lemma~\ref{lem:U=W+K} with $\kappa=\kappa_0$ gives the condition~2,
and condition~3 is nothing but Lemma~\ref{lem:positivity} with $s=0$ and $t=kT$.
We can thus infer from Corollary~\ref{cor:abstract} the existence of a unique triplet $(\Lambda_0,f_0,\phi_0)\in \R\times L^1_+(\R^d)\times L^\infty_+(\R^d)$ such that
   \[Sf_0=\Lambda_0f_0,\quad S^*\phi_0=\Lambda_0\phi_0,\quad \text{with}\quad \|\phi_0\|_{L^\infty}=\langle\phi_0,f_0\rangle=1.\]
We set $\lambda_F=\dfrac{\log\Lambda_0}{T}$ and, for $t\in[0,T]$,
\[f_t=e^{-\lambda_Ft}S_{t,0}f_0\qquad\text{and}\qquad \phi_t=e^{-\lambda_F(T-t)}S^*_{t,T}\phi_0,\]
so that we have $f_T=f_0$ and $\phi_T=\phi_0$, and we then define $f_t$ and $\phi_t$ for $t\in\R$ by $T$-periodicity.
These families verify the conditions $\langle\phi_t,f_t\rangle=\|\phi_0\|=1$ for all $t\in\R$ and the fact that $t\mapsto e^{\lambda_Ft}f_t$ satisfies~\eqref{eq:main} and $t\mapsto e^{-\lambda_Ft}\phi_t$ satisfies~\eqref{eq:main-dual}.
Besides, the convergence~\eqref{eq:discrete-ergo} in Corollary~\ref{cor:abstract} also reads
\[\big\|e^{-n\lambda_F T}S_{nT,0}f-\langle\phi_0,f\rangle f_0\big\|_{L^1}\leq C \zeta^n\left\|f-\langle\phi_0,f\rangle f_0\right\|_{L^1}.\]
Using~\eqref{eq:expo-growth} and the semiflow property, we classically extend this inequality to the continuous setting to get the existence of $C>0$ and $\rho>0$ such that
\begin{equation}\label{eq:conv-expo}
{\big\|e^{-\lambda_F(t-s)}S_{t,s}f-\langle\phi_s,f\rangle f_t\big\|}_{L^1}\leq Ce^{-\rho (t-s)}{\|f-\langle\phi_s,f\rangle f_s\|}_{L^1}
\end{equation}
for all $f\in L^1(\R^d)$ and all $t\geq s\geq0$.
We have proved the first part of Theorem~\ref{thm:main}, namely the existence of the Floquet eigenelements and the non-constructive exponential convergence.

\

Now we assume that {\bf(H$\v$)}-{\bf(H$a$)}-{\bf(H$q+$)} are satisfied and we aim at quantifying the constants $C$ and~$\rho$ by applying Theorem~\ref{thm:Harris} to $S$.
We start with some quantified estimates on $\lambda_F$ and $\phi_0$.
Using~\eqref{eq:expo-growth} and Lemma~\ref{lem:g_0} we have the quantitative estimate
\[\frac{\log\kappa_0}{T}\leq\lambda_F\leq \hat q+\esssup a.\]
Applying~\eqref{eq:HarrisR1R2} to $\phi=\phi_0$ we get the lower bound
\begin{equation}\label{eq:phi0-1}
\phi_0=e^{-\lambda_F T}S^*\phi_0\geq c\,e^{-(\hat q+\esssup a) T}\langle\phi_0,\1_{B(0,R_1)}\rangle \1_{B(0,R_2)}
\end{equation}
for some constructive constant $c=c(R_1,R_2)>0$.
Next, applying~\eqref{eq:LinfL1} to $\phi=\phi_0$ and $\lambda=\Lambda_0$ yields, since $(W^*-\Lambda_0)^{-1}K^*\phi_0=\phi_0$,
\[\langle\phi_0,h\rangle\geq(1-\alpha){\|\phi_0\|}_{L^\infty}=1-\alpha.\]
Due to the expression of $h$ in~\eqref{eq:h}, we infer that there exists $R_h>0$ such that
\begin{equation}\label{eq:phi0-2}
\langle\phi_0,\1_{B(0,R)}\rangle>0,\qquad\forall R>R_h.
\end{equation}
Combining~\eqref{eq:phi0-1} and~\eqref{eq:phi0-2}, we deduce that for any $R>0$ there exists $c_R>0$ such that
\begin{equation}\label{eq:phi0-3}
\phi_0\geq c_R\1_{B(0,R)}.
\end{equation}
We also provide an upper bound on $\phi_0$.
From the Duhamel formula~\eqref{eq:Duhamel1_dual} applied to $\phi_0$, and using~\eqref{eq:expo-growth}, we have
\[\phi_0(x)=e^{-\lambda_F T}S^*\phi_0(x)\leq \Phi(x):=e^{-\lambda_FT+\int_0^T\overline a(X_{s,0}(x))ds}+e^{(\hat q+\esssup a-\lambda_F)T}\hat q\int_0^T e^{\int_0^s\overline a(X_{s',0}(x))ds'}ds.\]
By invoking the dominated convergence theorem, we see that this function $\Phi$ tends to $0$ as $|x|\to\infty$.

\medskip

Now we prove that $\widetilde S=e^{-\lambda_FT}S$ satisfies the Lyapunov condition~\eqref{eq:Lyapunov} by computing for any $(t,x)\in\R^{1+d}$
\[\cL^*_t\1(x)=a(t,x)+\int_{\R^d}q(t,x,y)dy\leq\overline a(x)+\hat q.\]
Since $\lim_{|x|\to\infty}\overline a(x)=-\infty$, we can find $R,\theta>0$ such that
\[\cL^*_t\1\leq -1-|\lambda_F|+\theta\1_{B(0,R)}\leq-1-|\lambda_F|+\frac{\theta}{c_R}\phi_0 .\]
Injecting this inequality in~\eqref{eq:weak-mild} and using Grönwall's inequality, we obtain that
\[{\|\widetilde Sf\|}_{L^1}=e^{-\lambda_FT}{\|Sf\|}_{L^1}\leq e^{-T}{\|f\|}_{L^1}+\Theta\langle\phi_0,|f|\rangle\]
for any $f\in L^1$ and some constructive $\Theta>0$.
Now that the Lyapunov condition~\eqref{eq:Lyapunov} is proved, we look at Harris's condition~\eqref{eq:Harris}.
By duality, the result in Corollary~\ref{cor:Harris} is equivalent to the fact that for any $f\in L^1_+$
\begin{equation}\label{eq:HarrisR1R2-dual}
Sf\geq c\langle\1_{B(0,R_2)},f\rangle\1_{B(0,R_1)}.
\end{equation}
Consider $A>0$ such that $A>\Theta/(1-e^{-T})$.
 For $f\in L^1_+$ such that $\|f\|_{L^1}\leq A\langle\phi_0,f\rangle$ and $R>0$ we write
\begin{align*}
\langle\phi_0,f\rangle=\int_{B(0,R)}\phi_0f+\int_{B^c(0,R)}\phi_0f&\leq\langle\1_{B(0,R)},f\rangle+\Big(\sup_{B^c(0,R)}\Phi\Big){\|f\|}_{L^1}\\
&\leq\langle\1_{B(0,R)},f\rangle+\Big(\sup_{B^c(0,R)}\Phi\Big)A\langle\phi_0,f\rangle.
\end{align*}
Choosing $R$ large enough so that $A\sup_{B^c(0,R)}\Phi\leq1/2$, we obtain
\[\langle\1_{B(0,R)},f\rangle\geq\frac12\langle\phi_0,f\rangle\]
and~\eqref{eq:HarrisR1R2-dual} with $R_2=R$ and $R_1=A$ then gives
\[\widetilde S f\geq \langle\phi_0,f\rangle g_A\qquad\text{with}\quad g_A=\frac{1}{2}c\,e^{-\lambda_FT}\1_{B(0,A)}.\]
Since~\eqref{eq:phi0-3} implies that $\langle\phi_0,g_A\rangle$ is bounded from below by $\frac{1}{2}c_Ac|B(0,A)|e^{-\lambda_FT}$, Theorem~\ref{thm:Harris} applies to $S$ and guarantees that the constants $C$ and $\rho$ in~\eqref{eq:conv-expo} can be quantified.
This finishes the proof of Theorem~\ref{thm:main}. 

\

\paragraph{\bf Acknowledgments.}

The authors have been supported by the ANR project NOLO (ANR-20-CE40-0015), funded by the French Ministry of Research.

This work was partially funded by the European Union (ERC, SINGER, 101054787). Views and opinions expressed are however those of the author(s) only and do not necessarily reflect those of the European Union or the European Research Council. Neither the European Union nor the granting authority can be held responsible for them.


\end{document}